\documentclass{article}

\usepackage{graphicx}
\usepackage{amssymb}
\usepackage{amsmath, amsthm}
\usepackage[affil-it]{authblk}
\usepackage{color}

\usepackage{hyperref}

\usepackage{subcaption}

\usepackage[T1]{fontenc}

\usepackage{amsmath,amssymb}

 \usepackage{amsthm}
 
  \newtheorem{theorem}{Theorem}
  \newtheorem{remark}[theorem]{Remark}
     \newtheorem{definition}[theorem]{Definition}

\begin{document}

\title{Building up a model family for inflammations}

\author[1,$\ddagger$]{Cordula Reisch}

\author[2]{Sandra Nickel}
            
\author[2]{Hans-Michael Tautenhahn}

\affil[1]{Technische Universit\"at Braunschweig, Universit\"atsplatz 2, 38106 Braunschweig, Germany}
            
\affil[2]{Clinic for Visceral, Transplantation, Thoracic and Vascular Surgery; Leipzig University Hospital, Liebigstrasse 20, 04103 Leipzig, Germany, sandra.nickel@medizin.uni-leipzig.de, hans-michael.tautenhahn@medizin.uni-leipzig.de}

\affil[$\ddagger$]{Corresponding author: c.reisch@tu-braunschweig.de}

\date{December 08, 2023}

\maketitle

\begin{abstract}
The paper presents an approach for overcoming modeling problems of typical life science applications with partly unknown mechanisms and lacking quantitative data: 
A model family of reaction diffusion equations is built up on a mesoscopic scale and uses classes of feasible functions for reaction and taxis terms. 
The classes are found by translating biological knowledge into mathematical conditions 
and the analysis of the models further constrains the classes. 
Numerical simulations allow comparing single models out of the model family with available qualitative information on the solutions from observations. 
The method provides insight into a hierarchical order of the mechanisms. 
The method is applied to the clinics for liver inflammation such as metabolic dysfunction-associated steatohepatitis (MASH) or viral hepatitis where reasons for the chronification of disease are still unclear and time- and space-dependent data is unavailable. 
\end{abstract}



\section{Introduction}

Modeling life science processes with only partly known mechanisms and actors is challenging. 
The challenge becomes even harder if there is only few data or only qualitative data available for validating the models. 
However, mathematical models for not fully understood complex problems provide insight into the relevance of mechanisms and allow testing hypotheses in-silico. 

In this paper, a family of models for the dynamics of liver inflammation is formulated, analyzed, simulated, and considered in a clinical context. 
Liver inflammation like steatohepatitis or viral hepatitis often evolves into chronic inflammations which may lead to severe secondary diseases such as cirrhosis.

The mechanisms for the evolution of the inflammation course are not fully understood and therefore treatment is less effective.
Mathematical models provide a scientific environment for testing hypotheses on the importance of involved cells and mechanisms. 
A hierarchical order of the mechanisms is a priori unavailable because the mechanisms are coupled and only in a few cases isolatable for single experiments. 
Models on a mesoscopic scale allow an abstraction from the unknown processes on the cell scale and allow comparing the modeling results to qualitative data like pathological images of the liver.

The model family describes the interactions and propagation of different influencing variables on inflammation. In particular, the following factors are currently considered in the model: macrophages (copper cells), viruses, CD4$^{+}$ T cells, cytotoxic T cells, and cytokines. 
Different approaches for known mechanisms are used and models of different complexity are analyzed and compared. 
As the focus lies on the model family, not every model will be discussed in all detail and not all mechanisms will be used in the models. 
The results for the model family provide insight into the hierarchical order of the used mechanisms. 

\subsection{Literature overview}

The literature covers three topics: First, medical and biological observations on liver inflammations are given. 
Details on the known mechanism and involved cells follow in Sec.~\ref{sec:biological}. 
Second, models for inflammations are presented with a highlight on models for viral liver infections leading to inflammation. 
Third, some results on reaction diffusion equations are given. 

The leading causes of liver inflammation worldwide are viruses, metabolic diseases, autoimmune diseases, and alcohol. If left untreated over a long period of time, all of these factors can lead to cirrhosis of the liver, which is in fact life-threatening, \cite{kazankov_role_2019}.

The course of viral liver infections begins with the infection itself, followed by a symptom-free incubation period.
During an active phase, the immune system fights against the virus in the liver tissue. 
After the active phase, healing infection courses differ from chronic. 
While the amount of virus and the activity of the immune system reduces to zero in a healing infection course, the virus persists in the liver and the immune system remains active during a chronic infection course, \cite{schiff_schiffs_2018}.
Inflammation is the process of an active immune system reacting to infected tissue. 
If the inflammation continues even for reduced viral loads, the inflammation is called chronic. 

Inflammation occurs in different tissues caused by virus, bacteria, or autoimmune diseases. 
This paper focuses on inflammation caused by viral liver infections like hepatitis B. 
The reasons for chronic infections are still partly unknown, compare \cite{bowen_adaptive_2005,thomas_experimental_2016}, and data is rare. 
Qualitative data is given by pathological images of removed tissue, \cite{kanel_pathology_2017}.
The pathological images show a partial inhomogeneous spread of T~cells in the regarded tissue parts.  
The amount of T cells is higher next to portal fields.  
During acute hepatitis B infections, up to 95\% of the liver cells are infected with the virus, \cite{kanel_pathology_2017, schiff_schiffs_2018}.
During chronic phases around one-third of the liver cells are infected,  \cite{kanel_pathology_2017, schiff_schiffs_2018}.
These observations will be used as qualitative data in Sec.~\ref{sec:simulations} for evaluating the modeling results.  

Mathematical models for inflammations are available on different length scales and for different inflammation types. 
In \cite{aston_new_2018} models using ordinary differential equations are compared. 
The models describe the total amount of free virus, healthy, and infected liver cells in the whole liver by compartment models. 
Variations change the reaction functions or take delay into account, compare \cite{nangue_analysis_2022}. 
The same components are regarded in \cite{rezounenko_viral_2018} but the diffusion of all three cell types are added. 
In this approach, even the liver cells are diffusing, which is not realistic. 
A more realistic model that is adapted to hepatitis B infections is presented in \cite{tadmon_global_2021}.
The modeled cells are next to three types of liver cells, free virus and B cells and only the free virus is diffusing. 
This leads to a coupled ODE-PDE system for which stationary states are analyzed and chronic infections are regarded. 
The presented models describe the dynamics of the liver cells and the free virus. 
As available data is given by pathological images showing the T~cells, the results are only partly comparable with the space-dependent images. 

Mathematical models need to be spatially resolved for comparison with pathological images. 
Additionally, the T cells need to be modeled as these are the only cell types next to liver cells displayed in the pathological images. 
Reaction diffusion equations are one possibility for modeling interacting cells in space and time. 
They allow a mathematical analysis and in some cases a prediction of the solution behavior. 

Some results on reaction diffusion equations are summarized here. 
Reaction diffusion equations are used for modeling many different applications like morphogenesis \cite{turing_chemical_1952}, the spread of populations \cite{murray_mathematical_2002} or chemical reactions \cite{schnakenberg_simple_1979}. 
Depending on the reaction functions and the diffusion parameters, the models show different types of solutions like traveling waves, (Turing) patterns, blow-ups, or leveling behavior, \cite{perthame_parabolic_2015}.
In some cases, these solutions can be proven and predicted analytically.  
In the light of chronic and healing infection courses of liver inflammations, two types of solutions are particularly interesting:
Solutions tending towards zero describe healing infection courses and chronic infection courses are mathematically described by solutions with a tendency to stationary and spatially inhomogeneous distributions. 
As there are some results on leveling solutions, c.f. \cite{smoller_shock_1994}, the contradiction of these results can be used for gaining requirements for solutions not tending towards a leveled state, compare \cite{reisch_entropy_2020}.

\subsection{Overview of this paper}

The model family is formulated in Sec.~\ref{sec:modellfamily}. 
Biological information is translated into classes of feasible reaction functions (Sec.~\ref{subs:reaction}) and of taxis terms (Sec.~\ref{subsec:taxis}).
Subsections \ref{subs:ex_reaction} and \ref{subs:ex_taxis} give examples of these classes.
The model family is analyzed in Sec.~\ref{sec:modelfamily} and requirements on the classes are formulated in Sec.~\ref{subsec:requirements}.
Sec.~\ref{sec:simulations} presents numerical simulations and compares them with qualitative observations. 
A conclusion follows in Sec.~\ref{sec:conclusions}.

\section{Model family}
\label{sec:modellfamily}

A model family consists of models describing an application using the same modeling approach but differing in the number or complexity of included mechanisms.
First, the model type of reaction diffusion equations and boundary conditions are discussed. 
Afterward, biological information on the dynamics of liver infections is given and is translated into conditions on function classes.

\subsection{Reaction Diffusion Equations}\label{subs:reacdiff}

The models have the general form 
\begin{align}\label{eq:reacdiff_general}
\mathbf{q}_{,t} = \mathbf{F}(t, \mathbf{x}, \mathbf{q}) + \nabla \cdot \mathbf{D}(\mathbf{x}, \mathbf{q}, \nabla \mathbf{q}) ,
\end{align}
where $t>0$, $\mathbf{x} \in \Omega \subset \mathbb{R}^d$ and $\Omega$ is a model for the regarded part of the liver. 
Following, the dimension of space is $d \in \{ 1, 2, 3\}$, depending on the chosen simplification.

The vector-valued function $\mathbf{q}: [0, \infty) \times \Omega \rightarrow \mathbb{R}^n$ gives the time and space-dependent amount of the acting cells and virus. 
The function $\mathbf{F}$ describes the reactions between the cells and the virus, as well as the growth or decay of any described substances. 
The term $\mathbf{D}(\mathbf{x}, \mathbf{q}, \nabla \mathbf{q})$ models any diffusive or chemotactic effects caused by gradients of the substances or outer influence. 

The reaction diffusion in Eq.~(\ref{eq:reacdiff_general}) is completed to an initial value problem by 
$\mathbf{q} (0, \mathbf{x})= \mathbf{q}_0 (\mathbf{x})$ for $\mathbf{x} \in \Omega$
with $\mathbf{q}_0 (\mathbf{x}) \geq \mathbf{0}$ for all $\mathbf{x} \in \Omega$. 
Additionally, we assume zero-flux boundary conditions
\begin{align}\label{eq:bc}
\frac{\partial \mathbf{q}}{\partial \mathbf{n}} & = \mathbf{0} & \text{ for }  \mathbf{x} \in \partial \Omega . 
\end{align}
The homogeneous Neumann boundary conditions are an approximation: 
The liver consists of eight different segments which are separated by infeasible tissue, \cite{schiff_schiffs_2018}.
There, zero flux conditions are reasonable because the exchange with the exterior of one segment happens via the blood vessels. 

If $\Omega$ is a part of such a liver segment, homogeneous Neumann boundary conditions are less intuitive. 
Assuming that the liver tissue in $\Omega$ is homogeneous, periodic boundary conditions are a reasonable choice. 
Biologically, liver inflammations are stronger in some parts of the liver than in others. 
Therefore, the assumption of similar structures next to each other is not perfect.

 As an approximation for nearly constant amounts of acting cells close to the boundary $\partial \Omega$, zero flux conditions are acceptable and reflect the higher importance of local reactions compared to diffusion over the boundary. 

\subsection{Biological motivation}\label{sec:biological}

Viral hepatitis such as hepatitis B is named after the virus type causing an infection and leading to inflammation of the liver tissue. 
Inflammation is a consequence of the reaction of the immune system to the infection. 

Hepatitis B infections show two typical infection courses.
In healing infection courses, the virus and the inflammation vanish after an acute inflammation with a strong active phase of the immune system.
Chronic infection courses have an active phase followed by a chronic phase with mild symptoms. 
In the chronic phase, virus and inflammation persist in the liver tissue. 
Chronic liver infections often lead to severe secondary diseases like liver cancer or cirrhosis. 

A short summarization of the cellular immune response follows. 
After the infection, the virus spreads out in the body and settles in the liver. 
There, it replicates and spreads out. 
During this phase, the adaptive immune reaction starts working. 
Dentritic cells process and opsonize virus particles to the T cells. 
The T~cells specify and develop to CD4$^+$ or CD8$^+$ T~cells. 
Some of the CD4$^+$ T~cells become T helper cells specializing in different types of T helper cells. 
One type of T helper cell enhances the cytokine concentration, which are signals caused by the interaction of tissue cells and virus.
Cytokines lead the way for T~cells towards the virus. 

The CD8$^+$ T~cells named cytotoxic T cells triggers the cell death of infected tissue cells.
Therewith, the tissue cell and the contained virus are destroyed. 
Inflammation is essential to neutralize the damaging stimulus and create the conditions for repair processes.
The inflammation is caused only indirectly by the virus. 

Next, the description of infection is split up into single mechanisms. 
From a model-theoretic point of view, this step is again a modeling process. 
We start with naming the actors of the dynamics.
Let $q_1= q_1(t, \mathbf{x})$ be the space and time dependent amount of the virus. 
According, $q_2=q_2(t, \mathbf{x})$ is the amount of T~cells, which can be divided into T helper cells $T_h=T_h(t, \mathbf{x})$ and cytotoxic T~cells $T_c=T_c(t, \mathbf{x})$.
Cytokines $q_3=q_3(t, \mathbf{x})$ are the third quantity.

The effects on each other are visualized by arrows or dependencies ($\sim$).

\begin{enumerate}
	\item[(M1)] virus replicate and spread out: $q_1 \rightarrow q_1 \nearrow$
	\item[(M2)] immune system reacts on virus by producing T~cells: $q_1 \rightarrow q_2 \nearrow$
	\item[(M3)] T helper cells produce and enhance cytokines: $T_h \rightarrow q_3 \nearrow$
	\item[(M4)] cytokines lead cytotoxic T~cells: $q_3 \sim T_c$
	\item[(M5)] cytotoxic T~cells kill virus: $T_c \rightarrow q_1 \searrow$
	\item[(M6)] in absence of virus, T~cells reduce (as well natural decay): $\neg q_1 \rightarrow q_2 \searrow$
	\item[(M7)] random spreading of T helper cells 
\end{enumerate}

\begin{figure}[htbp]
  \centering
\includegraphics[width=5.5cm]{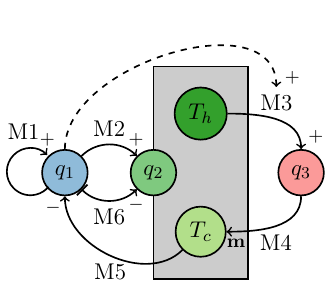}
  \caption{Mechanisms during an inflammation. 
  The variables are virus $q_1$, cytokines $q_3$ and T~cells $q_2$, divided into T helper cells $T_h$ and cytotoxic T~cells $T_c$. 
  Positive effects are displayed as $\footnotesize{ +}$, negative effects as $\footnotesize{-}$, and movement effects as \footnotesize{m}. 
  The spreading of the virus and T~cells is not visualized. 
  The dashed arrow shows the enhancement of mechanism M3 depending on $q_1$.}
   \label{fig:1}
\end{figure}

Fig.~\ref{fig:1} visualizes the mechanisms and their effect. 
There is a large loop effect: Mechanism M2 enhances the production of both T~cell types according to the virus amount $q_1$. 
The T helper cells $T_h$ produce, again depending on $q_1$, cytokines.
The cytokines $q_3$ control the cytotoxic T~cells $T_c$, which reduce the virus $q_1$.

Next,  the mechanisms are translated into assumptions on the reaction terms.

\subsubsection{Reactions}\label{subs:reaction}

The next step in building up a model family is the translation of the biological observations into mathematics. 
Fig.~\ref{fig:1} is an intermediate step for this process. 
There, the mechanisms are converted from a verbal description to a logical expression. 
Next, the logical expressions are translated into mathematical formulas, regarding more details of the verbal expression than in Fig.~\ref{fig:1}. 

Therefore, we specify the reaction diffusion system in \eqref{eq:reacdiff_general}.
$F_1$ is the first entry of the vector $\mathbf{F}$ and belongs to the partial differential equation for the virus $q_1$. 
Analogously, $F_2$ and $F_3$ belong to $q_2$ and $q_3$.  
If the T~cell population in the model is split up into the two sub-types $T_h$ and $T_c$, there are reaction functions $F_{T_h}$ and $F_{T_c}$.

Every reaction function $F_i$ consists of different single mechanisms $f_{i}^{Mj}$ splitting up the biological mechanisms occurring during an inflammation. 
The reaction function for the change of the virus is given by
\begin{align}\label{eq:F1}
F_1(t, \mathbf{x}, \mathbf{q}) = f_1^{M1}(t,\mathbf{x},\mathbf{q}) + f_1^{M5}(t, \mathbf{x}, \mathbf{q}),
\end{align}
where the first reaction function $ f_1^{M1}$ describes the growth and $ f_1^{M5}$ describes the decay of virus. 
The spreading of the virus, which is also part of mechanism M1, is described by a taxis term, which will be discussed in Sec.~\ref{subsec:taxis}.

If the T~cells are only regarded as general T~cells fulfilling all mentioned functions of T helper cells and cytotoxic T~cells, then the biological mechanisms affecting $q_2$ are
\begin{align*}
F_2(t, \mathbf{x}, \mathbf{q}) &= f_{2}^{M2}(t, \mathbf{x}, \mathbf{q})  + f_{2}^{M4}(t, \mathbf{x}, \mathbf{q}) + f_{2}^{M6}(t, \mathbf{x}, \mathbf{q}) .
\end{align*}
If the T~cell subtypes are split up into T helper cells $T_h$ and cytotoxic T~cells $T_c$, then the mechanisms affect directly the subtypes instead of the general T~cells, so 
\begin{align*}
F_{T_h}(t, \mathbf{x}, \mathbf{q}) & = f_{T_h}^{M2}(t, \mathbf{x}, \mathbf{q}) + f_{T_h}^{M6}(t, \mathbf{x}, \mathbf{q}), \\
F_{T_c}(t, \mathbf{x}, \mathbf{q}) & = f_{T_c}^{M2}(t, \mathbf{x}, \mathbf{q}) +f_{T_c}^{M4}(t, \mathbf{x}, \mathbf{q}) + f_{T_c}^{M6}(t, \mathbf{x}, \mathbf{q})  .
\end{align*}
So far, the reaction function for the cytokines consists of one mechanism,
\begin{align*}
F_3(t, \mathbf{x}, \mathbf{q}) &= f_{3}^{M3}(t, \mathbf{x}, \mathbf{q})  
\end{align*}
describing the production of cytokines by T helper cells in the presence of the virus. 

The mechanisms translate into conditions on the partial derivatives of the functions $f_i^{Mj}$.
By considering first the translation into partial derivatives and not directly into certain functions, the classes of reaction functions include more variety. 
The deductive modeling process from biological mechanisms to a mathematical model family is highlighted by these finer steps. 

In \cite{reisch_chemotactic_2019} the conditions
\begin{align*}
\frac{\partial F_1}{\partial q_2} <0, \quad \quad \frac{\partial F_2}{ \partial q_1} >0, \quad \quad \frac{\partial F_2}{\partial q_2}<0 
\end{align*}
are presented, which are now regarded in more detail. 
We start with the mechanisms for the virus. 
The growth of virus (M1) can be modeled by a function fulfilling
\begin{align}\label{eq:cond_m1}
\frac{\partial f_1^{M1}(t,\mathbf{x},\mathbf{q})}{\partial q_1} \geq 0.
\end{align} 
The equality with zero can be interpreted as a stopped growth of the virus, for example, due to an eliminated virus or due to saturation. 
The second mechanism affecting the virus is a reduction in dependency of cytotoxic T~cells. 
This mechanism can be expressed by the conditions on the partial derivatives
\begin{align}\label{eq:cond_m5}
\frac{\partial f_1^{M5}(t,\mathbf{x},\mathbf{q})}{\partial T_c} \leq 0 \quad \quad \mathrm{or} \quad \quad \frac{\partial f_1^{M5}(t,\mathbf{x},\mathbf{q})}{\partial q_2} \leq 0,
\end{align} 
depending on the chosen variables. 

Analogously, we find for the T~cells
\begin{align}\label{eq:partmech_q2}
\frac{\partial f_2^{M2}(t,\mathbf{x},\mathbf{q})}{\partial q_1} \geq 0  , \quad \quad \frac{\partial f_2^{M4}(t,\mathbf{x},\mathbf{q})}{\partial q_3} \geq 0 ,\quad \quad \frac{\partial f_2^{M6}(t,\mathbf{x},\mathbf{q})}{\partial q_2} \leq 0.
\end{align}
Mechanism M2 shows the effect that the immune system has on the virus by the production of T~cells. 
The more virus there is, the more T~cells are produced. 
The mechanism M4 shows that the amount of T~cells increases depending on the cytokines. 
In mechanism M6, the decay of T~cells is modeled. 

If the model considers a splitting up of the general T~cells into T helper cells and cytotoxic T~cells, then there are conditions on the mechanism functions for the special T~cells. 
For the T helper cells, the mechanisms M2 and M6 are relevant and the conditions in \eqref{eq:partmech_q2} hold for $f_{T_h}^{M2}$ and $f_{T_h}^{M6}$. 
The mechanisms M2, M4 and M6 are relevant for the cytotoxic T~cells and \eqref{eq:partmech_q2} hold for $f_{T_c}^{M2}$, $f_{T_h}^{M4}$ and $f_{T_c}^{M6}$.

The conditions for the reaction functions of the cytokines are
\begin{align} \label{eq:cond_3}
\frac{\partial f_3^{M3} (t, \mathbf{x}, \mathbf{q})}{ \partial T_h} \geq 0 \quad \quad \text{and} \quad \quad \frac{\partial f_3^{M3} (t, \mathbf{x}, \mathbf{q})}{ \partial q_1} \geq 0,
\end{align}
where the partial derivative with respect to the T helper cells is replaced by the partial derivative with respect to the T~cells $q_2$ in the case of not splitting up into subtypes. 

For all quantities, a natural decay will be discussed in Sec.~\ref{subs:ex_reaction}. 
The natural decay is a result of different biological processes and a small transport. 
For the virus, the natural decay can be interpreted as a reduced growth factor. 
Additionally, some processes have a natural saturation, for example, due to restricted resources.

The general conditions on the partial derivatives of the mechanism functions are one step of the modeling process providing the substance for building a model family.  
In the next section, conditions on the taxis terms are presented. 
Afterwards, particular reaction functions and taxis terms fulfilling the conditions are given and the effects of natural decay and saturation are discussed in more detail.

\subsubsection{Taxis}\label{subsec:taxis}

The taxis are not fully covered by the description in Fig.~\ref{fig:1} and therefore here explained in more detail.
The regarded domain $\Omega$ is a part of the liver, not containing, for example, the lymph.
Consequently, the details of the transport of T~cells into the liver is not part of the model. 
The domain $\Omega$ contains a subdomain $\Theta$ where the inflow of T~cells takes place. 
If $\Omega$ is a two-dimensional domain, $\Theta$ can be seen as a cut through a portal field of the liver structure. 

The boundary conditions prevent an exchange of virus or T~cells with the exterior. 
In this section, the taxis of any cells in the interior of $\Omega$ are described. 

The virus $q_1$ spreads out in the liver. 
The spreading is caused by two different mechanisms, \cite{goyal_modelling_2016}. 
Cell-to-cell transmission leads to a spreading with a small medium step size and the virus in one infected liver cell passes on to a neighbor cell after replication. 
The second mechanism is the diffusion of the virus through the extracellular space. 
The medium step size is typically larger than for the cell-to-cell transmission. 
The diffusion function in~\eqref{eq:reacdiff_general} for the $q_1$-component therefore consists of two terms
\begin{align*}
\mathbf{D}_1(\mathbf{x}, \mathbf{q},\nabla \mathbf{q}) = \mathbf{D}_1^\mathrm{ctc}(\mathbf{x}, \mathbf{q},\nabla \mathbf{q}) + \mathbf{D}_1^\mathrm{ecs}(\mathbf{x}, \mathbf{q},\nabla \mathbf{q}).
\end{align*}

The cell-to-cell transmission can be interpreted as a diffusion process with a small diffusion parameter.
This mechanism has the form 
\begin{align}\label{eq:D1_ctc}
 \mathbf{D}_1^\mathrm{ctc}(\mathbf{x}, \mathbf{q},\nabla \mathbf{q})  = \mathbf{D}_1^\mathrm{ctc} (\nabla q_1).
\end{align}

The spreading through the extracellular space might depend on the blood stream through the extracellular space and therefore on the direction and velocity of the blood. 
The results term would include this extracellular spread additionally to virus amount, for example
\begin{align}\label{eq:D1_flux}
 \mathbf{D}_1^\mathrm{ecs}(\mathbf{x}, \mathbf{q},\nabla \mathbf{q}) = \mathbf{D}_1^\mathrm{ecs}(\mathbf{x},q_1, \nabla q_1) . 
\end{align}

The cytotoxic T~cells  $T_c$ follow a cytokine gradient
\begin{align*}
 \mathbf{D}_{T_c}^\mathrm{chem} (\mathbf{x}, \mathbf{q},\nabla \mathbf{q})=\mathbf{D}_{T_c}^\mathrm{chem} ( T_c, \nabla q_3)
 \end{align*}
and the taxis depend on the cytotoxic T~cells and the cytokine gradient. 
 
 The T helper cells are spreading out by diffusion, so  $\mathbf{D}_{T_h}^\mathrm{diff} (\mathbf{x}, \mathbf{q},\nabla \mathbf{q})=\mathbf{D}_{T_h}^\mathrm{diff} ( T_h).$

For the general T~cells without a division into subtypes, both effects combine and the dependencies on $T_h$ and $T_c$ are replaced by a dependency on $q_2$.

The cytokines are signal cells without self-induced movement. 
Due to natural processes, a small diffusion is possible. 
Additionally, the transport of the cytokines by the bloodstream through the extracellular tissue might happen.
As a sum of those two mechanisms, the term for the cytokine $q_3$ is
\begin{align*}
\mathbf{D}_{3}(\mathbf{x}, \mathbf{q},\nabla \mathbf{q})= \mathbf{D}_3^\mathrm{diff}(\nabla q_3) + \mathbf{D}_3^\mathrm{ecs}(\mathbf{x}, q_3, \nabla q_3).
\end{align*}

We give examples for taxis terms in Sec.~\ref{subs:ex_taxis}.

\subsection{Examples of feasible reaction functions}\label{subs:ex_reaction}

The model in Eq.~\eqref{eq:reacdiff_general} is specialized for particular models of the model family. 
The functions are motivated by biology, see Sec.~\ref{sec:biological}, and still variable enough to cover various scenarios. 

Mechanism M1 describes the virus replication, and functions fulfilling Eq.~\eqref{eq:cond_m1} are 
\begin{align}\label{eq:mech1}
\begin{aligned}
f_1^{M1}(q_1)&= a_1 q_1 && (\text{unbounded}), \\
f_1^{M1}(q_1)&= a_1 q_1 (C_1 -q_1) && (\text{bounded by } C_1), \\
f_1^{M1}(q_1)&= a_1 q_1 (C_1 -q_1) \frac{q_1 - \varepsilon}{q_1 + \kappa} && (\text{bounded by }C_1 \text{ and Allee effect}).
\end{aligned}
\end{align}
The three examples are extensions of each other. 
The first function describes an unbounded exponential growth of the virus. 
The second function is bounded from above by a constant capacity $C_1$. 
The resulting growth is a logistic growth, starting exponentially and reaching saturation which describes that all liver cells are infected. 

The third example includes a strong Allee effect, \cite{allee_principles_1949}, modeling a small vanishing virus amount by a negative function for $q_1 < \varepsilon$. 
This effect can be interpreted as a small immune reaction acting on very small viral loads. 
A discussion of the difference due to the Allee effect can be found in \cite{reisch_chemotactic_2019}. 
The function with the Allee effect fulfills the condition in Eq.~\eqref{eq:cond_m1} only for values $q_1 \geq \varepsilon$.

The mechanism M2 describes the production of T~cells, both T helper cells and cytotoxic T~cells, if the virus is present. 
Possible mechanisms are 
\begin{align}\label{eq:mech2}
\begin{aligned}
f_2^{M2}(q_1)&= a_2 q_1 && (\text{local, unbounded}), \\
f_2^{M2}( q_1, q_2)&= a_2 q_1 (C_2 -q_2) && (\text{local, bounded by } C_2), \\
f_2^{M2}(\mathbf{x}, q_1)&= a_2 \chi_{\Theta}(\mathbf{x})  \int_\Omega q_1 \, \mathrm{d} \mathbf{x} && ( \text{global}), \\
f_2^{M2}(\mathbf{x}, q_1)&=a_2 \chi_{\Theta}(\mathbf{x}) (C_2 -q_2)   \int_\Omega q_1 \, \mathrm{d} \mathbf{x}  && ( \text{global, with saturation}). 
\end{aligned}
\end{align}
The first example depends only on the local virus amount. 
T~cells appear where the virus is. 
This effect is extended by an upper bound $C_2$ for the T~cells in the second example.
The local mechanism neglects the finer liver structure with portal fields. 
T~cells are produced outside of $\Omega$ and are transported via blood vessels into the liver. 
In the small-scale liver structure, the influx of T~cells takes place in portal fields $\Theta \subset \Omega$.

The subdomains $\Theta$ are part of the third and fourth examples for mechanism functions. 
The functions are explicitly space dependent via a function $\chi_{\Theta} (\mathbf{x})$ fulfilling
\begin{align}\label{eq:chi1}
\int_\Omega \chi_{\Theta} (\mathbf{x}) =1  \quad \text{ and } \quad \chi_{\Theta} (\mathbf{x})  \begin{cases} \geq 0 &\text{ for } \mathbf{x} \in \Theta,  \\ =0 &\text{ for } \mathbf{x} \in \Omega \setminus \Theta, \end{cases}
\end{align}
compare \cite{reisch_chemotactic_2019}.
The function $\chi$ can be a characteristic function for the subdomain $\Theta$, or, if required, a smooth function with positive values and a local support in $\Theta$.

The integral $\int_\Omega q_1 \, \mathrm{d} \mathbf{x}$ models a collection of the total viral load in the regarded part $\Omega$. 
The production of T~cells is proportional to the total viral load in examples three and four in Eq.~\eqref{eq:mech2}. 

If the T~cells are distinguished into T helper cells $T_h $ and cytotoxic T~cells $T_c$,
the proposed mechanisms are adapted by replacing $q_2$ by one of the subtypes with individual capacities $C_{T_h}$, $C_{T_c}$.
The requirement in Eq.~\eqref{eq:partmech_q2} is fulfilled for all examples. 

Mechanism M3 describes the production of cytokines by the virus and the effect of T helper cells increasing the cytokines. 
Examples of reaction functions are
\begin{align}\label{eq:mech3}
\begin{aligned}
f_3^{M3}(q_1)&= a_3 q_1  && (\text{only dependent on virus}), \\
f_3^{M3}( T_h)&= a_3 T_h && (\text{only dependent on T helper cells}), \\
f_3^{M3}( q_1, T_h)&= a_3 q_1 T_h && ( \text{unbounded}), \\
f_3^{M3}( T_h, q_3)&= a_3 T_h (C_3 - q_3) && ( \text{dependent on $T_h$ cells, bounded}), \\
f_3^{M3}( q_1, T_h, q_3)&= a_3 T_h q_1 (C_3 - q_3) && ( \text{bounded}), 
\end{aligned}
\end{align}
where $T_h$ can be replaced by $q_2$ if there is no division into subtypes of T~cells. 
 The examples differ in their dependencies on the involved cells and in the boundedness of cytokines at a certain place. 
 Dependency only on the virus neglects the mechanism that T helper cells increase the cytokines while functions independent of the virus neglect the mechanism that cytokines are produced by infected liver cells. 
  
 The conditions in Eq.~\eqref{eq:cond_3} are fulfilled. 
 The last two examples in Eq.~\eqref{eq:cond_3} depend on $q_3$ itself and give an upper bound. 
 By only regarding mechanism M3, the constant $C_3$ gives an upper bound. 
Further discussion on the boundedness of solutions follows in Sec.~\ref{sec:analytic}.
 
The examples for mechanism M3 in Eq.~\eqref{eq:mech3} can be hierarchically ranked. 
The more sub-mechanisms are included, the higher the ranking is. 
The bottom level contains only the first two examples.
The last example contains the most mechanisms and is therefore at the top of the ranking. 
Containing the most mechanisms is not always a feature for quality as the mechanism and their importance are unknown.

M4 describes the control of the movement of the cytotoxic T~cells, c.f. Sec.~\ref{subs:ex_taxis}.

The mechanism M5 describes the decay of the virus in dependence on cytotoxic T~cells. 
Possible functions for this mechanism are
\begin{align}\label{eq:mech5}
\begin{aligned}
f_1^{M5}(T_c)&= - a_5  T_c  && (\text{only dependent on cytotoxic T~cells}), \\
f_1^{M5}( q_1,T_c)&= -a_5 q_1 T_c && (\text{dependent on virus and cytotoxic T~cells}).
\end{aligned}
\end{align}
The condition in Eq.~\eqref{eq:cond_m5} is fulfilled by both functions in~\eqref{eq:mech5}. 

Again, the dependency on the cytotoxic T~cells can be replaced by a dependency on the general T~cells. 
The first example in Eq.~\eqref{eq:mech5} does not preserve the non-negativity of $q_1$. 
The analysis of the reaction functions follows in Sec.~\ref{sec:analytic} in more detail.

The mechanism M6 describes that the T~cell amount does not decrease if the virus is still present. 
It includes different effects: First, there is a decay of T~cells. 
Second, this decay is smaller, if there is more virus. 
The technical description of 'a lot of virus' depends on the boundedness of the virus $q_1$ by $C_1$. 
Examples are
\begin{align}\label{eq:mech6}
\begin{aligned}
f_2^{M6}(q_2)&= - a_6  q_2  && (\text{natural decay}), \\
f_2^{M6}( q_1,q_2)&= -a_6 q_2 (C_1 - q_1) && (\text{decay in absence of virus}), 
\end{aligned}
\end{align}
where the first example is included in the second because 
$$ f_2^{M6}( q_1,q_2)= -a_6 q_2 (C_1 - q_1) = -a_6 C_1 q_2 + a_6 q_1 q_2 = C_1 f_2^{M6} + a_6 q_1 q_2 $$
yields. 
Both examples in Eq.~\eqref{eq:mech6} fulfill the condition in Eq.~\eqref{eq:partmech_q2} on the partial derivative with respect to the T~cells. 
The second example in Eq.~\eqref{eq:mech6} consists of a natural decay and an increase depending on the virus and the T~cells. 
The increase can be interpreted as a benefit of the T~cells like in a predator-prey system or a host-pathogen system. 
Further details on this comparison are discussed in \cite{reisch_modeling_2019}. 
If the T~cells are divided into T helper cells and cytotoxic T~cells, mechanism M6 works in the same way by replacing $q_2$ by $T_h$ or $T_c$ respectively. 

Examples of the reaction functions of the different mechanisms are given in Eqn.~\eqref{eq:mech1}, \eqref{eq:mech2}, \eqref{eq:mech3}, \eqref{eq:mech5} and \eqref{eq:mech6}. 
The taxis describing mechanism M4 will be discussed in Sec.~\ref{subs:ex_taxis}.
Additionally, there is a natural decay for all quantities. 

The virus is decaying in the presence of (cytotoxic) T~cells, see mechanism M5 in Eq.~\eqref{eq:mech5}. 
Additionally, there is a small natural decay, which is already regarded in the growth, covered by mechanism M1 in Eq.~\eqref{eq:mech1}.
The growth including an Allee effect has a negative value for a very small amount of virus. 
This negative growth can be interpreted as a natural decay of very small amounts of virus. 
Besides, a natural decay can be expressed by a smaller parameter $a_1$. 
Consequently, the natural decay for the virus is already included in the example for mechanism M1 with the Allee effect. 
The other examples in Eq.~\eqref{eq:mech1} do not include the natural decay, interpreted as a local immune reaction, for small amounts of virus in the growth function. 
This results in different analytical properties, see Sec.~\ref{sec:analytic}. 

The reaction function of mechanism M6 for the T~cells contains a natural decay. 
In the case of a division of the T~cells into the subtypes of T helper cells and cytotoxic T~cells, the mechanism M6 is used in the same way for both subtypes. 

So far, a natural decay for $q_3$ like 
\begin{align}\label{eq:3_nd1}
f_{3}^\mathrm{nd} (q_3) =- a_\mathrm{nd} q_3,
\end{align}
is not contained in the mechanisms.
A possible constant decay is
\begin{align}\label{eq:3_nd2}
f_{3}^\mathrm{nd}  =\begin{cases} - a_\mathrm{nd} & \text{ if } q_3>0 ,\\
					0 & \text{ if } q_3=0 ,\end{cases} 
\end{align}
resulting in an unsteady reaction function. 

This section presents examples of reaction functions fulfilling the requirements from Sec.~\ref{subs:reaction} and describing mathematically the biologically motivated mechanisms.  
Next, we give examples of the movements of the different cells. 

\subsection{Examples of feasible taxis terms}\label{subs:ex_taxis}

The general reaction diffusion equation in Eq.~\eqref{eq:reacdiff_general} has a taxis term of the form $\nabla \cdot \mathbf{D} (t, \mathbf{x}, \mathbf{q}, \nabla \mathbf{q})$. 
This term was divided into terms for every component in Sec.~\ref{subsec:taxis}. 
In this section, the taxis terms for every component are specified and examples are given. 

The virus $q_1$ spreads out using cell-to-cell transmission and diffusion through the extracellular space. 
Possible dependencies for the two spreading mechanisms are given in Eq.~\eqref{eq:D1_ctc} and \eqref{eq:D1_flux}. 
Cell-to-cell transmission is a pure diffusion process, so
\begin{align*}
\mathbf{D}_1^\mathrm{ctc}(\nabla q_1)&= d_1^\mathrm{ctc} \nabla q_1. 
\end{align*}

The spreading of the virus through the extracellular space can be modeled as
\begin{align*}
\mathbf{D}_1^\mathrm{ecs}( \nabla q_1)&= d_1^\mathrm{ecs} \nabla q_1 && (\text{homogeneous diffusion}) ,\\
\mathbf{D}_1^\mathrm{ecs}(\mathbf{x}, \nabla q_1)&= d_1^\mathrm{ecs} A (\mathbf{x}) \nabla q_1 && (\text{space dependent diffusion}),
\end{align*}
where $A(\mathbf{x}) \in \mathbb{R}^{d \times d}$ is a positive definite, space-dependent matrix modeling the flux directions in the extracellular space at every point $\mathbf{x}$. 

In more detailed modeling approaches, the mechanisms could be additionally time-dependent and display structural deformation on cell-scale caused by inflammation. 

The taxis of T~cells is type-dependent. 
Cytotoxic T~cells are attracted by cytokines, whereas the taxis of T helper cells is dominated by diffusion. 
We start with the diffusion of the T helper cells which is analogous to the diffusion of the virus
 \begin{align*}
\mathbf{D}_{T_h}^\mathrm{diff}(\nabla T_h)&= d_{T_h}^\mathrm{diff} \nabla T_h. 
\end{align*}

The motion of the cytotoxic T~cells is dominated by chemotactic effects, which means that the cytotoxic T~cells follow a cytokine gradient. 
Examples are 
 \begin{align}\label{eq:chem_Tc}
 \begin{aligned}
\mathbf{D}_{T_c}^\mathrm{chem}(  \nabla q_3)&=  d_{T_c}^\mathrm{chem} \nabla q_3,\\
 \mathbf{D}_{T_c}^\mathrm{chem}( T_c, \nabla q_3)&=d_{T_c}^\mathrm{chem}T_c  \nabla q_3,
 \end{aligned}
\end{align}
where the first example in Eq.~\eqref{eq:chem_Tc} does not depend on the number of cytotoxic T~cells and might lead therefore to a negative value $T_c$.

If the T~cells are not divided up into the two subtypes, there might evolve conflicts with the dynamic of T~cells and cytokines. 
These problems will be discussed later. 
An addition of the motion of cytotoxic T~cells and T helper cells leads to a taxis term
 \begin{align*}
\mathbf{D}_{2}( q_2, \nabla q_2, \nabla q_3)=d_{2}^\mathrm{diff} \nabla q_2  +   d_{2}^\mathrm{chem}q_2  \nabla q_3.
\end{align*}
The cytokines are spreading by some diffusion $\mathbf{D}_{3}(  \nabla q_3)=d_{3}^\mathrm{diff} \nabla q_3  .$

Next, different models are built by using the variety in the model family. 

\section{Analysis of the model family}\label{sec:modelfamily}

The reaction functions and taxis terms in Sec.~\ref{subs:ex_reaction} and \ref{subs:ex_taxis} are combined to models of the model family. 
Afterward, analytical results for the models are discussed. 
This analysis gives further requirements on the variety of possible reaction and taxis functions. 
Numerical simulations of acceptable models follow in Sec.~\ref{sec:simulations}. 

\subsection{Examples for models}\label{sec:ex_models}

The model family includes models for different sets of cell types. 
The largest set is the combination $(q_1, T_h, T_c, q_3)$. 
Another possible set of participating cells is $(q_1, q_2, q_3)$.
The smallest model with interactions consists of virus and T~cells, $(q_1, q_2)$ and neglects the cytokines. 
 For certain analysis steps, different combinations or one-component models are meaningful as well, see \cite{reisch_chemotactic_2019, reisch_entropy_2020}.

\begin{remark}[Model 1]\label{rem:mod1}
A large model in our model family for the four components virus $q_1$, T helper cells $T_h$, cytotoxic T~cells $T_c$, and cytokines $q_3$ is
\begin{align}
q_{1,t} &= a_1 q_1 (C_1 -q_1) \frac{q_1 - \varepsilon}{q_1 + \kappa} -a_5 q_1 T_c + \nabla \cdot  [d_1^\mathrm{ctc} \nabla q_1 +d_1^\mathrm{ecs} A (\mathbf{x}) \nabla q_1] , \label{eq:mod1_q1}\\
T_{h,t} & = a_{2,h} \chi_{\Theta}(\mathbf{x}) (C_{T_h} -T_h)   \int_\Omega q_1 \, \mathrm{d} \mathbf{x} -a_{6,h} T_h (C_1 - q_1) +  \nabla \cdot [d_{T_h}^\mathrm{diff} \nabla T_h]  , \label{eq:mod1_Th}\\
T_{c,t} & = a_{2,c} \chi_{\Theta}(\mathbf{x}) (C_{T_c} -T_c)   \int_\Omega q_1 \, \mathrm{d} \mathbf{x} -a_{6,c} T_c (C_1 - q_1)-  \nabla \cdot [d_{T_c}^\mathrm{chem}T_c  \nabla q_3] , \label{eq:mod1_Tc}\\
q_{3,t} & = a_3 T_h q_1  - a_\mathrm{nd} q_3 + \nabla \cdot d_{3}^\mathrm{diff} \nabla q_3 , \label{eq:mod1_q3}
\end{align}
for $\mathbf{x} \in \Omega$ and $t>0$. 
Let $\mathbf{q} = (q_1, T_h, T_c, q_3)^{\mathrm{T}}$. 
Boundary conditions are like in Eq.~\ref{eq:bc} zero flux conditions. 
Non-negative initial values complete the problem. 
All parameters have positive values. 
The matrix $A$ is positive definite and $\chi_{\Theta}$ fulfills Eq.~\eqref{eq:chi1}.  
\end{remark}

Model 1 in Rem.~\ref{rem:mod1} uses all mechanisms and for every mechanism the example with the highest complexity in comparison to the other examples. 
Of course, there are many more detailed functions and models with many more components thinkable.
Next, a model with three interacting cell types is presented, compare \cite{reisch_chemotactic_2019}. 
\begin{remark}[Model 2]\label{rem:mod2}
A second model for inflammation contains the three components virus $q_1$, T~cells $q_2$ and cytokines $q_3$. 
The dynamics are given by 
\begin{align}
q_{1,t} &= a_1 q_1 (C_1 -q_1) \frac{q_1 - \varepsilon}{q_1 + \kappa} -a_5 q_1 q_2 + \nabla \cdot [ d_1^\mathrm{ctc} \nabla q_1 ] , \label{eq:mod2_q1}\\
q_{2,t} & = a_{2} \chi_{\Theta}(\mathbf{x})   \int_\Omega q_1 \, \mathrm{d} \mathbf{x} -a_{6} q_2 (C_1 - q_1) +  \nabla  \cdot [ d_{2}^\mathrm{diff} \nabla q_2 - d_{2}^\mathrm{chem}q_2  \nabla q_3 ], \label{eq:mod2_q2}\\
q_{3,t} & = a_3 q_1  - a_\mathrm{nd} q_3 + \nabla \cdot [d_{3}^\mathrm{diff} \nabla q_3 ], \label{eq:mod2_q3}
\end{align}
again completed with zero flux boundary conditions Eq.~\eqref{eq:bc} and non-negative initial values. 
The parameters are positive and $\chi_{\Theta}$ fulfills Eq.~\eqref{eq:chi1}. 
\end{remark}

As a further reduction, the cytokines are not regarded. 
The right circle in Fig.~\ref{fig:1} containing the mechanisms of the cytokines is neglected. 
The resulting model only contains the dynamics of virus and T~cells, compare \cite{kerl_reaction_2012, reisch_chemotactic_2019}. 

\begin{remark} [Model 3] \label{rem:mod3}
A two-component model of virus $q_1$ and T~cells $q_2$ is 
\begin{align}
q_{1,t} &= a_1 q_1 (C_1 -q_1) \frac{q_1 - \varepsilon}{q_1 + \kappa} -a_5 q_1 q_2 + \nabla \cdot  [d_1^\mathrm{ctc} \nabla q_1 ], \label{eq:mod3_q1}\\
q_{2,t} & = a_{2} \chi_{\Theta}(\mathbf{x})   \int_\Omega q_1 \, \mathrm{d} \mathbf{x} -a_{6} q_2 (C_1 - q_1) +  \nabla  \cdot [d_{2}^\mathrm{diff} \nabla q_2] . \label{eq:mod3_q2}
\end{align}
Again, the parameters are positive, $\chi_{\Theta}$ fulfills Eq.~\eqref{eq:chi1}, the initial conditions are non-negative and there are zero flux boundary conditions. 
\end{remark}

Model 3 in Rem.~\ref{rem:mod3} uses different examples of reaction functions including fewer effects than model 1 in Rem.~\ref{rem:mod1}.
This is not only a consequence of the reduction of modeled cell types but also depending on the used reaction functions. 

It is not possible to decide a priori, that means before an analytical analysis and some simulations, for the best model. 
This is caused by the uncertainty in choosing mechanisms and concrete reaction and taxis functions.
Further, it is a consequence of the missing definition for a 'best' model. 
  
With these three different models in mind, we start the analysis of the models with as general reaction functions as possible. 

\subsection{Analytical results}\label{sec:analytic}

The analysis of the reaction diffusion equations covers some basic properties. 

\subsubsection{Non-negativity of solutions}\label{subsec:nonneg}

The first property which will be discussed is the non-negativity of the solutions. 
As the solutions are interpreted as an amount of a certain cell type, negative values are meaningless in the light of application. 

The reaction diffusion equations consist of some reaction functions and divergence terms. 
If the divergence term considers only diffusive effects, the regarded system quantity will not become negative due to the diffusive term. 
In those cases, it is sufficient to regard the reaction terms and check whether it allows a negative quantity. 
We start with the discussion of the reaction functions and come afterward to the taxis. 

The reaction function for the virus $q_1$ consists of two mechanism functions for M1 and M5, see Eq.~\eqref{eq:F1}.
The mechanism M1 depends on the virus, so $f_{1}^{M1}=f_{1}^{M1}(t, \mathbf{x}, q_1)$. 
The explicit time and space dependency is only theoretical and not discussed.  

The mechanism M5 depends on the amount of (cytotoxic) T~cells, $T_c$ or $q_2$, and the amount of virus.
The dependency reads $f_{1}^{M5}=f_{1}^{M5}(t, \mathbf{x}, q_1, T_c)$, where $T_c$ can be replaced by $q_2$. 
Together, the reaction function is called $F_1(t, \mathbf{x}, q_1, T_c)$, and it specifies Eq.~\eqref{eq:F1}.
A necessary condition for the non-negativity of $q_1$ is
\begin{align}\label{eq:pos_q1}
F_1(t,\mathbf{x}, 0, T_c) \geq 0.
\end{align}
A function modeling the growth of a virus should fulfill 
\begin{align}\label{eq:m1_nonneg}
f_{1}^{M1}(t, \mathbf{x}, 0)=0.
\end{align}
All examples in Eq.~\eqref{eq:mech1} fulfill Eq.~\eqref{eq:m1_nonneg}, so there is no increase of virus if no virus is present.

The mechanism M5 depends not only on the virus itself but also on the amount of (cytotoxic) T~cells. 
Additionally, the requirement Eq.~\ref{eq:cond_m5} is a non-positive partial derivative with respect to the T~cells.
Therefore, $f_{1}^{M5}$ has negative values, and should be zero if $q_1=0$. 
Otherwise, there would be an automated cell death without any virus. 
Again, the requirement in Eq.~\ref{eq:pos_q1} transfers to a condition on $f_1^{M5}(t, \mathbf{x}, q_1, q_2)$ as
\begin{align}\label{eq:m5_nonneg}
f_1^{M5}(t, \mathbf{x}, 0 , q_2)=0.
\end{align}

Regarding the examples of reaction functions for mechanism M5 in Eq.~\eqref{eq:mech5}, only the second example fulfills these requirements. 
The first example depends only on the T~cells and can lead to a negative amount of virus if $q_1=0$. 

Next, we regard the reaction functions of the T~cells $q_2$.
The discussion holds as well for the subtypes T helper cells $T_h$ and cytotoxic T~cells $T_c$. 
The reaction function consists of two mechanisms, M2 for the production and M6 for the decrease. 
The production depends on the virus and, if a saturation effect is considered, on the amount of T~cells.
So, $f_2^{M2}= f_2^{M2}(t,\mathbf{x}, q_1, q_2)$. 
The decay of T~cells depends as well on the virus and on the amount of T~cells. 
While the dependency on $q_1$ is a fundamental requirement in mechanism M2, the dependency on $q_1$ in M6 is optional, compare Eqn.~\eqref{eq:mech2} and \eqref{eq:mech6}.
In total, the reaction function for the T~cells and the subtypes is a function $F_2 (t,\mathbf{x}, q_1, q_2)$.
As the function of M2 is non-negative, the requirement
\begin{align}\label{eq:pos_q2}
F_2(t,\mathbf{x}, q_1, 0) \geq 0
\end{align}
transfers directly to a requirement on the reaction function of mechanism M6, so 
\begin{align}\label{eq:pos_q2_m6}
f_2^{M6}(t,\mathbf{x}, q_1, 0)\geq0.
\end{align}
All examples in \eqref{eq:mech2} and \eqref{eq:mech6} fulfill Eq.~\eqref{eq:pos_q2_m6} and are therefore, concerning the non-negativity of solutions, suitable functions. 
Of course, other functions are thinkable as long as the non-negativity requirements are fulfilled. 
 
The reaction function $F_3$ of the cytokines $q_3$ consists of functions for the mechanism M3 and on the natural decay, compare Eqn.~\eqref{eq:mech3}, \eqref{eq:3_nd1} and \eqref{eq:3_nd2}.
Mechanism M3 describes the increase of cytokines and can depend on the virus, the T (helper) cells, and the amount of cytokines, so $ f_3^{M3}= f_3^{M3}(t, \mathbf{x}, q_1, q_2, q_3)$. 
The natural decay depends only on the amount of cytokines, either directly like in Eq.~\eqref{eq:3_nd1}, or indirectly like in Eq.~\eqref{eq:3_nd2}.
In total, the reaction function $F_3$ of the cytokines depends on all cells, so the non-negativity requirement reads
\begin{align}\label{eq:pos_q3}
F_3(t,\mathbf{x}, q_1,q_2, 0) \geq 0.
\end{align}
The function $f_3^{M3}$ describes an increase and therefore is non-negative. 
A function modeling the natural decay of cytokines should fulfill the requirement 
\begin{align}\label{eq:pos_q3_nd}
f_3^\mathrm{nd}(t,\mathbf{x}, 0) = 0
\end{align}
to model the natural decay. 
Both examples in Eq.~\eqref{eq:3_nd1} and \eqref{eq:3_nd2} fulfill the condition. 

As seen, most of the conditions on the reaction functions for preserving non-negativity transfer directly into conditions on single reaction functions. 
Most of the introduced examples for reaction functions fulfill these conditions. 

As a next step, the taxis terms are considered. 
We start with the diffusive movement, described by a term $\nabla \cdot [ d_i^{\mathrm{diff}} \nabla q_i ]= d_i^\mathrm{diff} \Delta q_i$. 

The diffusive spreading is well known, for example from the heat equation. 
It can not lead to a negative amount of $q_i$ if $q_i$ is non-negative for all $\mathbf{x} \in \Omega$. 
This is a consequence of the leveling nature of diffusion and can be shown by using the fundamental solution of the heat equation after estimating the reaction function $f$ as a non-negative value, compare \cite{evans_partial_2010}. 
The same follows in case of a diffusion term $\nabla \cdot [ d_i^{\mathrm{diff}}  A(\mathbf{x} ) \nabla q_i ]$, where $A$ is for all $\mathbf{x} \in \Omega$ positive definite. 
Consequently, the non-negativity of the solutions for the cells with only diffusive spreading is assured. 
 
In the case of chemotactic effects, like for the cytotoxic T~cells $T_c$, we show the non-negativity as well. 
An example of a reaction diffusion equation for $T_c$ is given by 
\begin{align*}
T_{c,t} = F_{T_c} (t, \mathbf{x}, q_1, T_c, q_3) + \nabla \cdot [  d_{T_c}^{\mathrm{chem}} T_c \nabla q_3 ],
\end{align*}
analogously for $q_2$.
As the reaction function fulfills $ F_{T_c} (t, \mathbf{x}, q_1, 0, q_3) \geq 0$, the non-negativity of $T_c$ is proven if the divergence term does not lead to $T_c (t, \mathbf{x}) <0$ for any $t>0$ and any $\mathbf{x} \in \Omega$. 
The chemotactical term is the second example in Eq.~\eqref{eq:chem_Tc}.

Under the assumption of a continuous function, $T_c(t, \mathbf{x})$ is zero before it might become negative. 
If $T_c$ is negative at a time $t$ and a place $\mathbf{x}$, then the term for the chemotactic effects is negative as well, due to the linear dependency on $T_c$. 
So, $T_c$ has non-negative values if the taxis term depends explicitly on $T_c$, analogously for $q_2$.

The first example in Eq.~\eqref{eq:chem_Tc} does not depend on $T_c$. 
This term might lead to a negative value of $T_c$, for example, if $T_c(t, \mathbf{x}) \equiv0$ for a certain time $t$ and all $\mathbf{x} \in \Omega$. 
Then, depending on the distribution of $q_3$, a negative diffusion $\Delta q_3$ is possible for a point $\mathbf{x}$. 
The reaction function can be zero at this point, for example, if $q_1(t, \mathbf{x}) \equiv 0$ as well. 
It is required, that the taxis term depends on the T~cells
\begin{align*}
\mathbf{D}_{T_c}^\mathrm{chem} = \mathbf{D}_{T_c}^\mathrm{chem} (T_c, \nabla q_3).
\end{align*}

Altogether, systems fulfilling the requirements on the reaction functions and on the taxis terms are candidates for well-suited models w.r.t. non-negative solutions. 
Of course, the requirements are only necessary and not sufficient in all cases.

 \subsubsection{Boundedness of the solutions}
 
The non-negativity of solutions, which is a lower bound of the solutions, was discussed in Sec.~\ref{subsec:nonneg}.
Now, an upper bound for the solutions is discussed. 
The effort to show the existence of an upper bound of the solutions depends strongly on the chosen reaction functions. 
Of course, mechanisms with a negative impact on the variable do not increase an upper bound. 
If the mechanisms with a positive impact on the variable are already bounded, like some examples in Eqn.~\eqref{eq:mech1}, \eqref{eq:mech2} and \eqref{eq:mech3}, then the solution is bounded.  
If one mechanism is unbounded, proving an upper bound requires more effort. 

From a biological point of view, the existence of upper bounds is a desired property. 
As the space in the liver is finite, there cannot be infinitely many particles. 

We discuss the boundedness of $q_1$, $q_2$, and $q_3$ separately. 
The existence of upper bounds for $T_h$ and $T_c$ follows along the discussion of $q_2$. 
In some cases, an upper bound of one variable is already required for formulating some other mechanisms. 

Again, we start with the boundedness of the virus $q_1$. 
As the virus needs liver cells for reproduction, the reproduction is limited by the number of liver cells in the domain $\Omega$. 
Additionally, there can be free virus in the extracellular space, but again, this space is limited as well. 
Consequently, a growth function of mechanism M1 $f_1^{M1}(t, \mathbf{x},q_1)$  should fulfill the requirement that 
\begin{align}\label{eq:cap_1}
f_1^{M1}(t, \mathbf{x}, C_1) = 0
\end{align}
for a capacity $C_1>0$. 
Only the first example in Eq.~\eqref{eq:mech1} does not fulfill this requirement. 
This example models an unbounded growth, which is not biologically realistic. 

The mechanism M6 depends on the existence of an upper bound for the virus. 

Next, we regard the cytokines $q_3$. 
The reaction function consists of a growth term and a natural decay. 
The growth function $f_3^{M3}$ may depend on the virus, the T (helper) cells, and the cytokines themselves. 

We start discussing the cases of $f_3^{M3}= f_3^{M3}(q_1)$. 
If, as discussed before, the virus is bounded by $C_1$, then the influence of $q_1$ on the increase of $q_3$ is bounded as well. 
Together with the natural decay, there will be a value $q_3(t, \mathbf{x}) =C_3$ where the natural decay equalizes the production of cytokines. 

For functions $f_3^{M3}= f_3^{M3}(q_1, q_2, q_3)$, the influence of $q_3$ gives an upper bound through a bounded increase.
In the case of $f_3^{M3}= f_3^{M3}(q_1, q_3)$, this upper bound is not required due to the interplay with the natural decay. 

So far, we did not discuss the boundedness of $q_2$. 
If $q_2(t, \mathbf{x}) \leq C_2$ for all $t>0$ and all $\mathbf{x} \in \Omega$, then there is an upper bound for $q_3$ as well. 
This follows by analogous arguments as for the dependency of the increase only on the virus. 
If $q_2$ is unbounded, the increase of cytokines $q_3$ is unbounded as well.

There are two cases of increasing T~cells regarding the boundedness, compare Eq.~\eqref{eq:mech2}. 
Either, the mechanism for the increase already includes a point-wise upper bound $C_2$ or not. 
If the mechanism does not include such as boundedness, then the discussion needs more effort, compare \cite{reisch_longterm_2022}.  
We refer to this longer discussion of a maximal value and formulate as a requirement that there is a value $\tilde{C}_2$ such that the time derivative $q_{2,t} (t, \mathbf{x}) $ becomes negative for $q_2 (t, \mathbf{x}) > \tilde{C}_2$.

This requirement is for example fulfilled for the model in Rem.~\ref{rem:mod3}, see \cite{reisch_longterm_2022}. 

Altogether, for gaining a realistic model, all variables should be in $L^\infty (\Omega_T)$. 
If the boundedness of the virus is included in its growth function, the boundedness of the T~cells $q_2$ and the cytokines $q_3$ follow in some cases.

  \subsubsection{Longtime behavior}
 
As the solution $\mathbf{q}= \mathbf{q}(t, \mathbf{x})$ gives the amount of virus, T~cells and cytokines at a time $t$ and a space $\mathbf{x}$, the time evolution of the solution can be interpreted as an infection course. 
The question answered in this section is what types of solutions can be expected and whether these can be predicted a priori. 
Solutions of reaction diffusion equations show different behavior like a tendency towards steady states, traveling waves, or blow-ups, compare \cite{perthame_parabolic_2015}.

In the light of modeling inflammations, requirements for the boundedness of solutions were formulated. 
These requirements and a bounded domain $\Omega$ lead to the boundedness of solutions in any $L^p$-norm. 
Consequently, blow-ups are neither desirable for modeling inflammations nor occurring. 
Solutions of traveling waves are common for modeling atherogenesis, a certain form of inflammation, compare \cite{ibragimov_mathematical_2005,volpert_elliptic_2014}.
This type of solution is not suitable for modeling liver infections with a chronic infection course, due to the biological observations of more stationary inflammations, compare \cite{schiff_schiffs_2018}.
The observed spread of T~cells can be interpreted as a spatially inhomogeneous stationary solution. 

\begin{remark}\label{rem:interpretation_soluations}
We interpret solutions tending towards a stationary state which is spatially inhomogeneous as chronic liver infections, compare \cite{kerl_reaction_2012}.
Solutions tending towards zero are interpreted as healing infection courses.
\end{remark}

Suitable models for inflammations should show both relevant solution types, depending on the used parameters in the model or on the domain $\Omega$. 
Analytically, it is desirable to have a priori knowledge about the long-term behavior of the solution. 
As the tendency towards a spatially homogeneous steady state is a well-analyzed behavior of reaction diffusion equations, compare for example \cite{smoller_shock_1994}, there is some hope to predict those solutions. 
The theory of reaction diffusion equations with solutions tending towards a spatially inhomogeneous steady-state solution is still an open research field. 

We discuss two approaches for a priori statements on the occurrence of decaying solutions. 
The first approach is based on \cite[Thm. 14.17]{smoller_shock_1994} and was as well used in \cite{kerl_reaction_2012} in the context of a selected model for inflammations. 

The theorem in \cite{smoller_shock_1994} yields for reaction diffusion equations ${\mathbf{q}_{,t}= \mathbf{F}(\mathbf{q}) + D \Delta \mathbf{q}}.$
Models with chemotactic effects are not covered directly by this theorem. 
The statement of the theorem requires some definitions. 

\begin{definition}[compare \cite{smoller_shock_1994}]
Let $\lambda$ be the first non-zero eigenvalue to the eigenfunctions of the negative Laplacian on $\Omega$ with zero flux boundary conditions. 
Further, let $d$ be the smallest eigenvalue of the diffusion matrix $D$. 
Let $\Sigma$ be the invariant domain of the solutions $\mathbf{q}$, $M= \max_{\mathbf{q} \in \Sigma} \lVert \frac{\partial \mathbf{F}}{\partial \mathbf{q}} \rVert_{\mathbb{R}^{n \times n}}$ and $\sigma = \lambda d - M$.
\end{definition}

\begin{remark}[compare \cite{smoller_shock_1994}]\label{rem:smoller}
The reaction diffusion system $\mathbf{q}_{,t}= \mathbf{F}(\mathbf{q}) + D \Delta \mathbf{q}$ has an invariant domain of the solutions $\Sigma$, $\mathbf{F}$ is smooth and $D$ is positive definite. 
If $\sigma >0$ yields, then 
\begin{align*}
\lVert \mathbf{q}(t, \mathbf{x}) - \bar{\mathbf{q}}(t) \rVert_{L^\infty(\Omega)} \leq c \mathrm{e}^{- \sigma t}
\end{align*} 
where $\bar{\mathbf{q}}(t)$ is the solution of $\bar{\mathbf{q}}_{,t} = \mathbf{F}(\bar{\mathbf{q}}) + \mathbf{g}(t)$ with $\mathbf{g}$ decaying exponentially. 

Further, if $D$ is diagonal, then $\mathbf{q}$ tends in $L^\infty(\Omega)$ towards $\bar{\mathbf{q}}$. 
\end{remark}
Similar results with slightly different requirements can be found in \cite{perthame_parabolic_2015}.

Next, the requirements of Rem.~\ref{rem:smoller} are checked for the models for inflammations. 
As the structure of the reaction diffusion equations only considers diffusion, models with chemotactic effects are not covered by the theorem.
This is a problem for all models with cytokines.  
Models with mechanisms like in the first example of Eq.~\eqref{eq:chem_Tc} can be written in the form of Rem.~\ref{rem:smoller} but those do not fulfill the conditions for non-negative solutions.

Solutions interpreted as chronic infections may occur if $\sigma <0$. 
This is a necessary condition for the occurrence of solutions tending toward a spatially inhomogeneous steady-state distribution. 
In \cite{kerl_reaction_2012} this is checked for a model similar to the one in Rem.~\ref{rem:mod3}.
In \cite{reisch_chemotactic_2019} some simplified models with chemotactic effects are analyzed. 

For the wide class of possible models in this paper, the theorem of Rem.~\ref{rem:smoller} is not well suited and fits only for models without chemotactic effects. 

A different approach for predicting the longtime behavior are entropy methods, compare \cite{jungel_entropy_2016}. 
As those methods require the definition of an entropy functional depending on the reaction and diffusion terms, a general statement on the use of entropy functionals is not possible. 
In \cite{reisch_entropy_2020} a simplified model for liver inflammations consisting only of one variable was analyzed using entropy methods. 
The main analytical result for the reduced model is that if the decay terms of the virus overrule the growth of the virus, then only solutions interpreted as healing infection courses occur. 
This result is not directly transferable to the analysis of a two or more component model as the interplay between the different components dominates the system behavior. 

Altogether, only a few conclusions can be drawn from the existing theorems on the longtime behavior of very general reaction diffusion systems. 

The next section gives an overview of the requirements of the reaction and taxis functions. 
Additionally, the requirements are checked for the three presented models.

\subsubsection{Summarize of the requirements}\label{subsec:requirements}

The requirements are summarized according to the cell types and the introduced mechanisms. 
First, the requirements R.1 belonging to the virus are presented. 

\begin{remark}\label{rem:req:q1}
The dynamics of the virus are described by 
\begin{align}
q_{1,t} & = f_1^{M1} (t, \mathbf{x},q_1) + f_1^{M5}(t, \mathbf{x},q_1, q_2) + \nabla \cdot \mathbf{D}_1 (\mathbf{x}, \nabla q_1)
\end{align}
where $f_1^{M1}$, $f_1^{M5}$ and $\mathbf{D}_1$ fulfill the requirements
\begin{enumerate}
	\item[(R.1.1)] $ f_{1,q_1}^{M1} \geq 0$ for $q_1 > \varepsilon >0$, compare Eq.~\eqref{eq:cond_m1}, 
	\item[(R.1.2)]  $f_1^{M1} (t, \mathbf{x},0)=0$, compare Eq.~\eqref{eq:m1_nonneg},
	\item[(R.1.3)] $f_1^{M1} (t, \mathbf{x},C_1)=0$ for a capacity $C_1 > \varepsilon >0$, compare Eq.~\eqref{eq:cap_1}
	\item[(R.1.4)] $f_1^{M5}(t, \mathbf{x},0, q_2) \geq 0 $, compare Eq.~\eqref{eq:cond_m5},
	\item[(R.1.5)] $f_{1,q_2}^{M5} \leq 0$, compare Eq.~\eqref{eq:m1_nonneg}, 
	\item[(R.1.6)] if $\mathbf{D}$ has the form $A \nabla q_1$, then $A$ is positive definite.
\end{enumerate}
If the general T~cells are divided into T helper cells and cytotoxic T~cells, the dependency of $f_2^{M5}$ on $q_2$ is replaced by a dependency on the cytotoxic T~cells $T_c$.
\end{remark} 

\begin{remark}\label{rem:req:q2}
The dynamics of the T~cells are described by 
\begin{align}
q_{2,t} & = f_2^{M2} (t, \mathbf{x},q_1, q_2) + f_2^{M6}(t, \mathbf{x},q_1, q_2 ) + \nabla \cdot \mathbf{D}_2 (\mathbf{x}, q_2, \nabla q_2, \nabla q_3)
\end{align}
which can be adapted for only T helper cells ($T_h$ instead of 2 as index) or cytotoxic T cells ($T_c$ instead of 2 as index).
In any case, the functions $f_2^{M2}$, $f_2^{M6}$ and $\mathbf{D}_2$  fulfill the requirements
\begin{enumerate}
	\item[(R.2.1)] $ f_{2,q_1}^{M2} \geq 0$, compare Eq.~\eqref{eq:partmech_q2}, 
	\item[(R.2.2)]  $f_2^{M2} (t, \mathbf{x},q_1,q_2)\geq 0$,
	\item[(R.2.3)] $f_{2,q_2}^{M6} (t, \mathbf{x},q_1,q_2)\leq 0$, compare Eq.~\eqref{eq:partmech_q2}
	\item[(R.2.4)]  $f_2^{M6} (t, \mathbf{x},q_1,0)\geq 0$, compare Eq.~\eqref{eq:pos_q2_m6},
	\item[(R.2.5)] there exists a value $C_2$ with $q_{2,t} \leq 0$ for all $q_2 > C_2$, 
	\item[(R.2.6)] if $\mathbf{D}_2$ has a function $\mathbf{D}_2^\mathrm{ecs}(t,\mathbf{x}, \nabla q_2)$ with $A \nabla q_2$, then $A$ is positive definite,
	\item[(R.2.7)] if $\mathbf{D}_2$ has a function $\mathbf{D}_2^\mathrm{chem} (q_2, \nabla q_3)$, then the dependency on $q_2$ is obligatory. 
\end{enumerate}
\end{remark}

As a last quantity, the cytokines are regarded. 

\begin{remark}\label{rem:req:q3}
The dynamics of the cytokines $q_3$ are described by 
\begin{align}
q_{3,t} & = f_3^{M3} (t, \mathbf{x},q_1, q_2, q_3) + f_3^\mathrm{nd}(t, \mathbf{x},q_3 ) + \nabla \cdot \mathbf{D}_3 (\mathbf{x},\nabla q_3)
\end{align}
where $f_3^{M3}$, $f_3^\mathrm{nd}$ and $\mathbf{D}_3$ fulfill the requirements
\begin{enumerate}
	\item[(R.3.1)] $ f_{3,q_1}^{M3} \geq 0$ and $ f_{3,T_h}^{M3} \geq 0$, compare Eq.~\eqref{eq:cond_3}, 
	\item[(R.3.2)]  $f_3^{M3} (t, \mathbf{x},q_1,q_2,q_3)\geq 0$,
	\item[(R.3.3)]   $f_3^\mathrm{nd} (t, \mathbf{x},q_3)\leq 0$ for all $q_3 \geq 0$,
	\item[(R.3.4)]  $f_3^\mathrm{nd} (t, \mathbf{x},0)= 0$, compare Eq.~\eqref{eq:pos_q3_nd},
	\item[(R.3.5)] there exists a value $C_3$ with $q_{3,t} \leq 0$ for all $q_3 > C_3$, 
	\item[(R.3.6)] if $\mathbf{D}_3$ has the form $A \nabla q_3$, then $A$ is positive definite.
\end{enumerate}
If the T~cells are not split up, the partial derivative in (R.3.1) with respect to $T_h$ is replaced by the partial derivative with respect to $q_2$. 
\end{remark}

The requirements in Rem.~\ref{rem:req:q1}, \ref{rem:req:q2}, \ref{rem:req:q3} are reasonable conditions on the function of a thinkable model for inflammations. 
Of course, they are neither sufficient nor necessary for formulating a model describing inflammation in a general context. 

For the first application, the requirements are checked for the models 1,2 and 3. 

\begin{remark}\label{rem:checkm1}
Model 1 (Rem.~\ref{rem:mod1}) fulfills all requirements (R.1), (R.2), (R.3).
\end{remark}

\begin{remark}\label{rem:checkm2}
Model 2 (Rem.~\ref{rem:mod2}) fulfills all requirements (R.1), (R.2) and (R.3) under the assumption that $q_1 \leq C_1$.
Due to the requirements, the boundedness of the virus is preserved by the dynamics if the initial conditions are locally bounded by $C_1$. 
The boundedness of the T~cells $q_2$ is a consequence of the observations on the two-component model in \cite{reisch_longterm_2022}. 
\end{remark}

\begin{remark}\label{rem:checkm3}
Model 3 in Rem.~\ref{rem:mod3} fulfills the requirements (R.1) and (R.2).
The boundedness of the T~cells $q_2$ was shown in \cite{reisch_longterm_2022}. 
\end{remark}

In total, all three models fulfill the requirements, even though, in some cases, the proof of the boundedness is more complicated and a priori predictions are not always possible.

\subsection{Numerical results}\label{sec:simulations}

Simulations provide a visual understanding of models and possible solutions. 
In this section, simulations for the models 1-3 are presented. 

The simulations use a semi-discretization of the space coordinates for $\Omega = (0,1) \times (0,1)  \subset \mathbb{R}^2$. 
Afterward, standard solvers for ordinary differential equations are used. 
As the motion terms include Laplacians, the resulting ordinary differential equations of the space discretization may become stiff due to the growing quotient of the smallest and the largest eigenvalue of the discretized Laplacian. 
The number of discretization points is fixed by $\Delta x=0.05$ and, in case of long calculation time for solving the ordinary differential equations, a solver for stiff equations is used.

This section focuses on two points. 
First, the two solution types interpreted as chronic and as healing infection courses are presented. 
Second, there is a comparison of the three models using two to four different cell types. 
A leading question is whether even a model covering only two of four cell types can provide any insight. 

The initial conditions are chosen as $q_1(0, \mathbf{x}) \equiv 1$, $q_2(0, \mathbf{x}) = T_c(0, \mathbf{x})= T_h(0,\mathbf{x}) \equiv 0$ and $q_3(0, \mathbf{x}) \equiv 0.1$ in all simulations.  
The function $\chi_\Theta (\mathbf{x})$ is the normalized characteristic function for $\Theta$. 
Most of the parameters are fixed for all simulations, see Table~\ref{tab:para}.

\begin{table}[h!]
\caption{Parameter values for all three models which are chosen as constant if the mechanism is included. }
\begin{center}
\begin{tabular}{c c c c c c c c c  c c c c c c c c c c c }
 $a_1$ & $C_1$ & $\varepsilon$ & $\kappa$  & $d_1^\mathrm{ctc}$ & $a_{2,h}$ & $C_{T_h}$ & $a_{6}$ & $d_{T_h}^\mathrm{diff}$ & $C_{T_c}$ & $a_3$ &  $a_\mathrm{nd}$ & $d_3^\mathrm{diff}$\\[0.2em] \hline
 1 & 1 & 0.05& 0.01 & 0.6 &2 &8 & 0.2 & 0.9  &15 &0.8 &0.6 &0.5 
\end{tabular}
\end{center}
\label{tab:para}
\end{table}%

\begin{figure}[htb]
  \centering
  \includegraphics[width=\textwidth]{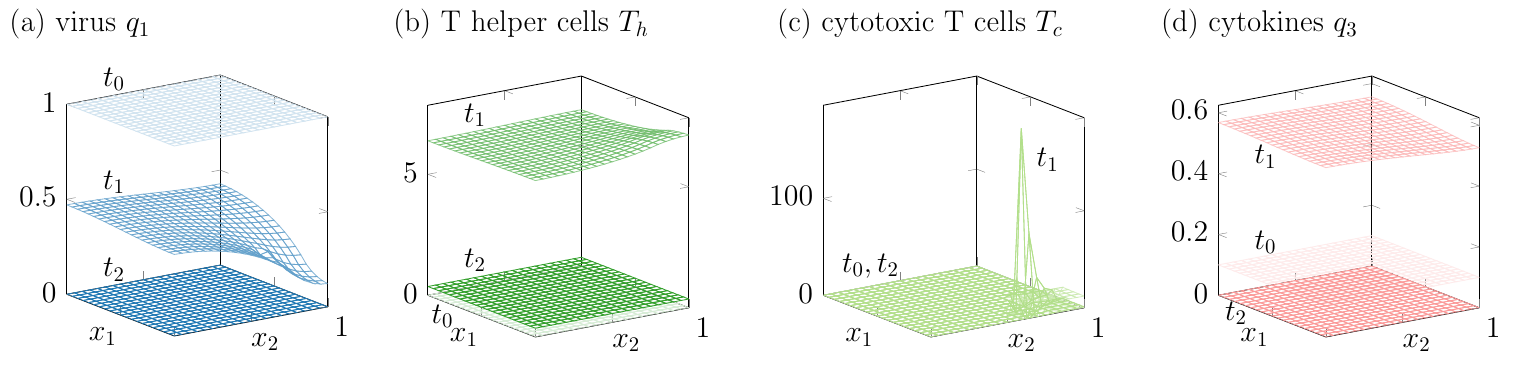}
    \includegraphics[width=\textwidth]{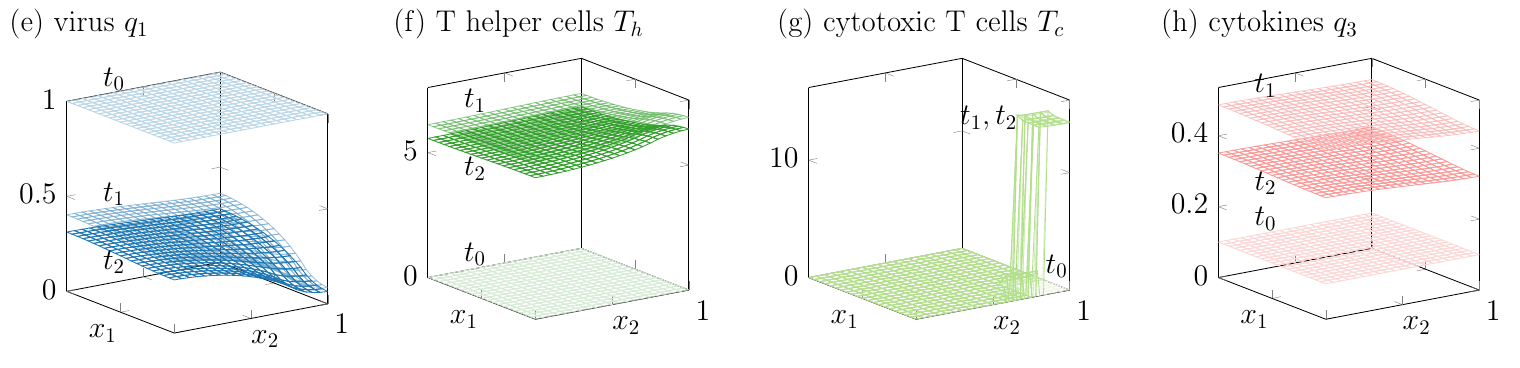}
  \caption{Model 1 shows solutions interpreted as healing infection courses (upper row) and as chronic infections (lower row). 
  Different time steps are in different shadings from $t_0=0$ in light over $t_1 =16$ to $t_2=40$ in dark.  }
    \label{fig:mod1healing}
\end{figure}

Model 1 shows solutions with a tendency towards zero and solutions with a tendency towards a spatial inhomogeneous steady state, depending on one parameter.
 In this case, the parameter $d_{T_c}^\mathrm{diff}$ changes the behavior, see Fig.~\ref{fig:mod1healing} and Table~\ref{tab:modchang} for the used parameter values.
The spread of the virus is modeled by diffusive cell-to-cell transmission.  
Extracellular movement through a blood flow is not regarded, $d_1^\mathrm{ecs}=0$. 

\begin{table}[htp]
\caption{Parameter values for the changed parameters in bold.
}
\begin{center}
\begin{tabular}{l|lll|lll|ll}
        & \multicolumn{3}{c|}{model 1}                 & \multicolumn{3}{c|}{model 2}                                 & \multicolumn{2}{c}{model 3}        \\
        & $a_5$ & $a_{2,c}$ & $d_{T_c}^\mathrm{chem}$  & $a_5$                      & $a_{2}$ & $d_{2}^\mathrm{chem}$ & $a_5$ & $a_{2}$                    \\[0.3em] \hline
healing & 2     & 2         & {\bfseries 1} & {\bfseries 1}   & 2       & 1                     & 0.5   & {\bfseries  2}   \\
chronic & 2     & 2         & {\bfseries 8} & {\bfseries 0.5} & 2       & 1                     & 0.5   & {\bfseries  0.7}
\end{tabular}
\end{center}
\label{tab:modchang}
\end{table}%

The simulations of model 1 show that there is a remarkable difference in the modeled spreading of T helper cells $T_h$ and cytotoxic T~cells $T_c$. 
While the spread of T helper cells is modeled by diffusion in model 1, compare Rem.~\ref{rem:mod1}, the motion of cytotoxic T~cells is modeled by chemotaxis. 
Consequently, the cytotoxic T~cells $T_c$ in Fig.~\ref{fig:mod1healing}(h) spread out less and the immune reaction is less effective in reducing the viral load. 
In the chronic case, the number of T helper cells $T_h$ reduces after an active phase but remains high for all time. 

\begin{figure}[htbp]
  \centering
\includegraphics[width=0.75\textwidth]{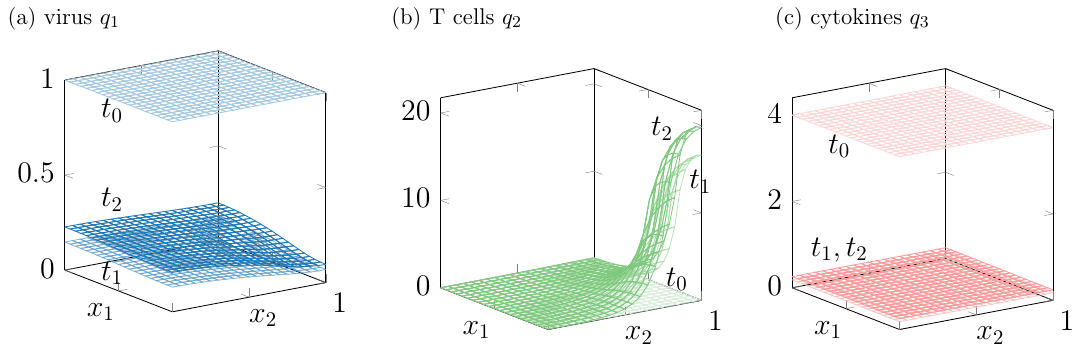}
  \caption{Model 2 shows solutions interpreted as chronic infection courses. 
  Different time steps are in different shadings from $t_0=0$ in light over $t_1 =32$ to $t_2=80$ in dark.
  The portal field $\Theta$ is in the right corner in (b). }
    \label{fig:mod2chronic}
\end{figure}

In model 2, there is only one type of T~cells spreading out by diffusion and following the gradient of the cytokines. 
The parameter under variation in this case is $a_5$, see Table~\ref{tab:modchang}, which describes the effectiveness of the T~cells. Fig.~\ref{fig:mod2chronic} shows a solution interpreted as chronic infection course. 
The virus remains in the whole liver with a spatial inhomogeneous spreading. 
The cytokines remain at a high level everywhere. 

\begin{figure}[htbp]
  \centering
\includegraphics[width=0.8\textwidth]{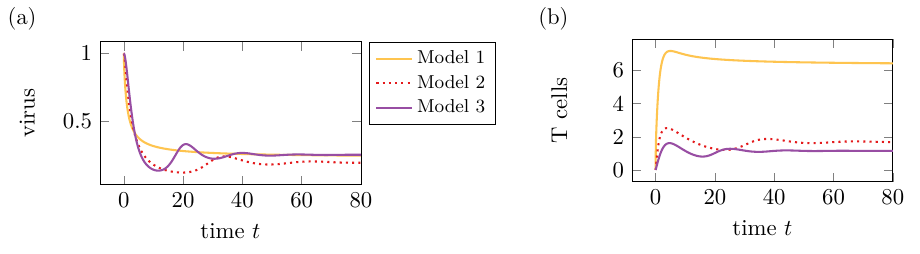}
  \caption{$L_1$ norms of the solutions of model 1, 2 and 3 for the parameter sets in Table~\ref{tab:modchang} interpreted as chronic infection course. 
  (a) total virus (b) total T~cells.
 }
   \label{fig:mod123}
\end{figure}

The changed parameter for model 3 is $a_2$, see Table~\ref{tab:modchang}, which regulates the inflow of T~cells through the portal field in dependency on the total virus at a time $t$. 
The results of the simulation of model 3 are comparable to the simulation results of models 1 and 2. 
A comparison of all three models for a chronic course is in Fig.~\ref{fig:mod123}. 
As the T~cells $q_2$ are a sum of the T helper cells $T_h$ and the cytotoxic T~cells $T_c$, Fig.~\ref{fig:mod123} compares the $L_1$ norms of $q_1$ in (b) of $q_2$ and $T_h+T_c$. 

The amount of virus tends in all three models to comparable values due to the small changes of parameters. 
The total amount of T~cells is comparable for models 2 and 3 but different for model 1 due to the differentiation between T helper cells and cytotoxic T~cells with different maximal values. 
The curves of models 2 and 3 show small oscillations in the first time but a damping towards a constant value later on. 

\subsection{Clinical impact}\label{sec:simulations}
The mathematical model for liver inflammation developed and described in this thesis may have significant clinical implications as it provides valuable insights into the underlying dynamics and mechanisms of liver inflammation caused by viruses. Some potential clinical improvements include a deeper understanding of the disease and its progression, supporting clinical decisions through early detection and diagnosis, and optimizing treatment. All of this could contribute to personalized diagnosis and treatment of viral hepatitis in the future.
More specifically, this mathematical model can help to detect subtle changes in influencing variables or physiological parameters associated with liver inflammation at an early stage. This early detection allows for rapid intervention and treatment that can potentially prevent the progression of liver inflammation to cirrhosis. Such predictive models can help physicians use their resources more efficiently by identifying patients who are likely to require more intensive monitoring or intervention. This is of paramount importance in healthcare systems with limited resources. Mathematical models provide a platform for simulating the dynamic processes involved in liver inflammation. Understanding disease progression at a mechanistic level can contribute to the development of new therapeutic targets and strategies or prevent inflammation from becoming chronic. More and more patients are demanding shared decision making, where the visualization and explanation derived from the mathematical model can be used to educate patients about their disease. Better patient understanding can improve treatment adherence and lifestyle change, which has a positive impact on long-term outcomes.
In summary, a mathematical model for liver inflammation can provide clinicians with valuable tools to improve diagnosis, treatment planning, and patient outcomes. Its integration into clinical practice has the potential to usher in a new era of precision medicine for liver disease.

\section{Conclusions}
\label{sec:conclusions}

A deductive modeling approach was presented for the application of liver infections leading to inflammations. 
This life science application has unknown mechanisms leading to chronic infection courses for which only qualitative data is available. 
Building up a model family allows to gain a deeper understanding of the involved mechanisms. 
The model family consists of reaction diffusion equations and the reaction or taxis functions can be chosen from a class of feasible functions. 
The feasible classes fulfill different biologically motivated properties. 
Analytical investigations restricted the function classes further.  

In the case of inflammation modeling, three models were chosen out of the model family and all three models reproduce both different infection courses depending on some parameter values. 
Consequently, already the smallest of the three models is able to cover basic observations even if less information is used than in the more complex models of the model family. 
Depending on the modeling purpose, a suitable model from the model family can be chosen. 

\section{Acknowledgments}

HM Tautenhahn disclosed receipt of the following financial support for the authorship of this article: This work was supported by:

Federal Ministry of Education and Research; ATLAS [031L0304C] HM Tautenhahn

German Research Foundation SIMLIVA [465194077] HM Tautenhahn

German Research Foundation [428832822] HM Tautenhahn\\

\noindent This research did not contain any studies involving animal or human participants, nor did it take place on any private or protected areas. No specific permissions were required for corresponding locations.

\linespread{1}

  \bibliographystyle{elsarticle-num} 
    \bibliography{references.bib}





\end{document}